\documentclass[a4paper,11pt]{article}
\usepackage{bbm}
\usepackage{amsfonts}
\usepackage{amsthm,amsmath}
\usepackage[integrals]{wasysym}
\usepackage{mathrsfs}
\usepackage{amsmath,multicol,amsopn}
\usepackage{latexsym,amssymb}
\usepackage{amsbsy, amstext}
\usepackage{graphicx}    
\usepackage{indentfirst} 
\usepackage{graphicx}
\usepackage{subfigure}
\usepackage{color}
\def\be{\begin{array}}
\def\en{\end{array}}

\numberwithin{equation}{section}

\newtheorem{theorem}{Theorem}[section] 
\newtheorem{definition}{Definition}[section]

\newtheorem{lemma}{Lemma}[section]
\newtheorem{remark}{Remark}[section]

\newcommand{\R}{{\mathbb R}}
\newcommand{\pozhehao}{\kern0.3ex\rule[0.8ex]{1.5em}{0.095ex}\kern0.3ex}

\begin{document}
\title{\LARGE\bf{Global Regularity to the liquid crystal flows of Q-tensor model}
}
\date{}
\author{Zhi Chen$^{1,2}$\thanks{{\small E-mail:  zhichenmath@ahnu.edu.cn (Z. Chen)}},
~~
 Elide Terraneo $^{2}$\thanks{{\small E-mail:elide.terraneo@unimi.it (E. Terraneo)}},~
 $^{1}$\thanks{{\small Corresponding author: E-mail: elide.terraneo@unimi.it }},\\
{\small 1. School of Mathematics and Statistics, Anhui Normal University,}\\
{\small Wuhu 241002 Anhui, P. R. China}\\
{\small 2. Dipartimento di Matematica
Universit\`{a} di Milano,}\\
{\small Via C. Saldini 50 20133 Milano, Italy}\\
}\maketitle

\begin{center}
\begin{minipage}{15.5cm}

{\bf Abstract.} In this paper we investigate a forced incompressible Navier-Stokes equation coupled with a parabolic type
equation of Q-tensors in a domain $U\subset\R^3.$ In the case $U$ is bounded, we prove the existence of a global strong solution when the initial data are sufficiently small, improving a result in Xiao's paper [J. Differ. Equations 2017].
The key tool of the proof is a {maximum principle.} Then, we establish also a result
of continuous dependence of solutions on the
 initial data. Finally, if  $U=\R^3,$ based on a result of Du, Hu and Wang [Arch. Rational Mech. Anal. 2020], we give an interesting regularity criterium  just via the $\dot{B}^{-1}_{\infty,\infty}$ norm of $u$ and the $L^\infty$ norm of the initial data $Q_0$.



\end{minipage}
\end{center}
\vskip 6mm

\section{Introduction and main results}
 \setcounter{equation}{0}
\par
In this paper, we study the model for a viscous incompressible liquid crystal flow proposed by Beris and Edwards
\cite{MR1352465}, which one can find in the physics literature, such as  in \cite{MR2084306,PhysRevE.67.051705}. An important feature of the
model is the molecular orientation
of the liquid crystals described by the so-called tensor $Q.$ More precisely,
let $U\subset\R^3$ be a bounded domain and $t>0,$ then for any  $(t,x)\in\R^+\times U,$ the tensor $Q(t,x)$ is a
symmetric traceless matrix, namely, for any $t>0$ and $x\in U$,  $Q(t,x)\in S^{(3)}_0 $, where
\begin{equation*}
   S^{(3)}_0=\{Q\in\R^{3\times3};\ Q=Q^T,\  trace[Q]=0\}.
\end{equation*}
Its motion is described through the velocity field $u=u(t,x),\ t\in R^+,\  x\in U\subset\R^3,$ by the coupled system:
\begin{equation}\label{gs11111}
\begin{cases}
\partial_tQ+(u\cdot\nabla)Q-S(\nabla u,Q)=\Gamma(\Delta Q-\mathcal{L}[\partial F(Q)]),\\
\partial_tu+(u\cdot\nabla)u-\nu\Delta u+\nabla P= \mathrm{div}(\tau+\sigma),\\
\mathrm{div}\, u=0,
\end{cases}
\end{equation}
with
\begin{equation*}
  \mathcal{L}[A]:=A-\frac{1}{3}trace[A]{Id}
\end{equation*}
denoting the projection onto the space of traceless matrices.
$P(t,x)$ is the scalar pressure, $\nu>0$ is the viscosity coefficient and $\Gamma>0$ is the elasticity relaxation constant. Moreover,
\begin{equation*}
  S(\nabla u,Q):=(\xi \epsilon(u)+w(u))\left(Q+\frac{1}{3}Id\right)+\left(Q+\frac{1}{3}Id\right)(\xi\epsilon(u)-w(u))-2\xi
  \left(Q+\frac{1}{3}Id\right)trace[Q\nabla u].
\end{equation*}
where $\epsilon(u):=\frac{1}{2}(\nabla u+\nabla^tu)$ and $w(u):=\frac{1}{2}(\nabla u-\nabla^tu)$
 are the symmetric
and antisymmetric part of the velocity gradient, respectively, and $\xi \in \R$.
The symmetric part of the additional stress tensor is given by:
\begin{equation}\label{tau}
\begin{aligned}
  \tau_{\alpha\beta}:=&-\xi\left(Q_{\alpha\gamma}+\frac{\delta_{\alpha\gamma}}{3}\right)H_{\gamma\beta}-\xi H_{\alpha\gamma}
  \left(Q_{\gamma\beta}+\frac{\delta_{\gamma\beta}}{3}\right)\\
  &+2\xi\left(Q_{\alpha\beta}+\frac{\delta_{\alpha\beta}}{3}\right)Q_{\gamma\delta}H_{\gamma\delta}-\partial_\beta Q_{\gamma\delta}\partial_\alpha Q_{\gamma\delta},
\end{aligned}
\end{equation}
and the antisymmetric part is
\begin{equation*}
  \sigma_{\alpha\beta}:=Q_{\alpha\gamma}H_{\gamma\beta}-H_{\alpha\gamma}Q_{\gamma\beta}.
\end{equation*}
where $H:=\Delta Q-\mathcal{L}[\partial F(Q)]$ and  we use the Einstein summation convention. Moreover,  the potential $F$ is a function defined on symmetric matrices that
will be introduced later and
the term $S(\nabla u,Q)$ describes how the flow gradient
rotates and stretches the order-parameter $Q$. The constant $\xi$ denotes the ratio of the tumbling and aligning effects on the
liquid crystals exerted by a shear flow.

\par
Throughout the paper, we denote  $A:B=A_{\alpha\beta}B_{\alpha\beta},$ $A\cdot B=trace(AB)$. Note that $A:B=A\cdot B$ for $A,B\in S^{(3)}_{0}$.  Moreover, we use the Frobenius norm of a matrix defined by
$|Q|:=\sqrt{trace(Q^2)}=\sqrt{Q_{\alpha\beta}Q_{\alpha\beta}}$ for $Q\in S^{(3)}_{0}$.
 Let $|\nabla Q|^2=\partial_\gamma Q_{\alpha\beta} \partial_\gamma Q_{\alpha\beta}$ and $|\Delta Q|^2:=\Delta Q_{\alpha\beta}\Delta Q_{\alpha\beta}.$
 We denote
 the norms of the usual Lebesgue spaces $L^p$ by $
\|u\|^p_{L^p}=\int_{U}|u|^p dx,~\hbox{for}~1\leq p<\infty$.
$C_i$, $C$  denote different
positive constants in different places.
 Finally, we note
$
  (u\otimes u)_{\alpha\beta}=u_\alpha u_\beta,$ and $(\nabla Q\otimes\nabla Q)_{\alpha\beta}=\partial_\alpha Q_{\gamma\delta}\partial_\beta Q_{\gamma\delta}.
$

\par
Several competing theories exist which attempt to capture the complexity of nematic liquid crystals. In 1960s,  Ericksen and his collaborators
proposed the so called Ericksen-Leslie equation for the underlying flow field which is given by the Navier-Stokes equation with an
additional forcing term  (see \cite{MR137403}).  This theory is very successful in many aspects and it matches well with the experimental observations.
However, this model cannot describe the biaxial nematic liquid crystals. So a new and more comprehensive theory is introduced by P. G. de Gennes \cite{de1993physics}, in which he uses $N\times N$ symmetric, traceless matrix Q as the new parameter, which is called $Q$ tensor.
\par
As mentioned, in this paper we study the model proposed by Beris and Edwards \cite{MR1352465}.
This  model has been {extensively and continuously studied} by many authors in recent years. In the case of polynomial potentials
of the form
\begin{equation}\label{CZ1}
  F(Q)=\frac{a}{2}|Q|^2-\frac{b}{3}trace[Q^3]+\frac{c}{4}|Q|^4, ~a,b,c\in\R,
\end{equation}
Paicu and Zarnescu \cite{MR2864407} established the existence of global weak solutions in dimension $N=2$ or $3$ with $\xi=0.$
Moreover, solutions with higher order regularity and the weak-strong uniqueness for $N=2$ is studied. Later, they extend these results
to the case when $|\xi|$ is sufficiently small in \cite{MR2837493}. In \cite{MR3485967}, with no additional smallness assumptin on $\xi$, they
establish the existence of global weak solution and long-time behavior of system with nonzero $\xi$ in the two-dimensional periodic domain.
The uniqueness of weak solutions to the Cauchy problem in $\R^2$ was proved in \cite{MR3599422,MR3576270}.
\par
On the other hand, concerning strong solutions, Abels, Dolzmann, and Liu \cite{MR3252809} proved the existence of a local strong solution and global weak solution
with higher regularity in time in the case of inhomogeneous mixed Dirichlet and Neumann boundary condition, in a bounded domain, without
smallness assumption on $\xi.$ Later, Liu and Wang in \cite{MR3927580} proved the existence and uniqueness of local in time strong solutions to the system with an anisotropic elastic energy. In the paper \cite{MR3912677}, Schonbek and Shibata obtain a global strong solution in $\R^N,(N\geq2)$ for  more general function $F.$ Moreover, they give an optimal decay rate for the $Q$ tensor in $L^p, ~~1\leq p\leq\infty.$
In the case of polynomial potentials like \eqref{CZ1}, Xiao \cite{MR3582195} proved the global well posedness for the simplified model with $\xi=0$ in a bounded domain under some assumptions on the coefficients. He constructed a solution in the maximal $L^p-L^q$ regularity class.
\par
Regarding  the decay of the solutions, it is worth pointing out that Dai et al. give a decay rate for the solutions $(u,Q)$
with general $F.$ Specifically,
in  \cite{MR3486134} they consider functions
\begin{equation*}
  F: \R^{3\times3}_{sym}\rightarrow(-\infty,\infty],
\end{equation*}
where $\R^{3\times3}_{sym}$ is the set of $3\times 3$ symmetric matrices. They suppose that  $F\in C^2(O),$  for some open set  $O\subset\R^{3\times3}_{sym}$
 containing the isotropic state $Q\equiv0.$
They suppose also that $Q=0$ is the (unique) global minimum of $F$ in $O,$ specifically,
\begin{equation}\label{CW1}
  F(0)=0,~~ F(Q)>0 ~~\text{for any}~~Q\in O\setminus\{0\}
\end{equation}
and
\begin{equation}\label{CW2}
  \partial F(Q):Q\geq 0 ~~\text{whenever}~~Q\in B_{r_1}~~\text{or}~~Q\in O\setminus B_{r_2}.
\end{equation}
where  $B_{r_1}, B_{r_2}$
with $r_1<r_2$ are two balls  such that
\begin{equation*}
Q=0\in B_{r_1}\equiv\{|Q|<r_1\}\subset\{|Q|\leq r_2\}\subset O.
\end{equation*}
\par

They show that the solutions of the system \eqref{gs11111} tend to the isotropic state at the rate $(1+t)^{-\frac{3}{2}}$ as $t\rightarrow\infty$ when $\xi=0$.
Moreover, if $F$ is as \eqref{CZ1}, they obtain the exponential decay.
It is worth pointing out  that they introduce an interesting maximum principle.  Note that any $F(Q):=\frac{a}{2}|Q|^2-\frac{b}{3}trace[Q^3]+\frac{c}{4}|Q|^4,$ with
 $a>0$, satisfy the conditions    \eqref{CW1}, \eqref{CW2},  at least  in a neighborhood of $0$.
\par
In the present paper, we consider $\xi=0$ which denotes  the  situation when the molecules only tumble in a shear flow but they are not aligned by the flow. Moreover, we consider the polynomial potential $F(Q)=\frac{a}{2}|Q|^2-\frac{b}{3}trace[Q^3]+\frac{c}{4}|Q|^4,$
as in the paper of Paicu and Zarnescu \cite{MR2864407}.
In this case, $\mathcal{L}[\partial F(Q)]=aQ-b\left[Q^2-\frac 13 Id \ trace (Q^2)\right]+cQ~ trace(Q^2)$, and 
 since $\xi=0$, the system \eqref{gs11111} reduces to the following:
\begin{equation}\label{cz22}
\begin{cases}
\partial_tQ+(u\cdot\nabla)Q-w(u) Q+Q w(u)=\Delta Q-aQ+b\left[Q^2-\frac 13 Id \ trace (Q^2)\right]-cQ~ trace(Q^2),\\
\partial_tu+(u\cdot\nabla)u=\Delta u-\nabla P+ \mathrm{div}(\nabla Q \otimes \nabla Q)+ \mathrm{div}(Q~ \Delta Q- \Delta Q~  Q),\\
\mathrm{div}\ u=0,
\end{cases}
\end{equation}
where,  without loss of generality, we set  $\nu=\Gamma=1$.
Moreover,  the following initial-boundary conditions are prescribed,
\begin{equation}\label{cz33}
  (u,Q)|_{t=0}=(u_0,Q_0),~~~(u,\partial_{\bf \nu} Q)|_{\partial U}=(0,\mathbb{O}),
\end{equation}
where ${\bf \nu}$ is the outward normal vector on $U$ and $ \mathbb{O}$ is the zero matrix.
\par
The main result of our paper is an improvement of the result of global existence of a strong solution
in Xiao's paper  \cite{MR3582195}.
In order to present Xiao's result, let us first introduce the definition of strong solution.
\begin{definition}\label{CH1}
For $T>0$ and $1<p,q<\infty,$ $M^{p,q}_T$ denotes the set of triplets $(u,Q,P)$ such that
\begin{equation*}
  u\in C([0,T],D^{1-\frac{1}{p},p}_{A_q})\cap L^p(0,T;W^{2,q}(U)\cap W^{1,q}_0(U)),
\end{equation*}
\begin{equation*}
  \partial_tu\in L^p(0,T;L^q(U)), \mathrm{div}u=0,
\end{equation*}
\begin{equation*}
  Q\in C([0,T];B^{3-\frac{2}{p}}_{q,p})\cap L^p(0,T;W^{3,q}(U)),~~\partial_tQ\in L^p(0,T;W^{1,q}(U))
\end{equation*}
\begin{equation*}
  P\in L^p(0,T;W^{1,q}(U)), \int_{U}Pdx=0.
\end{equation*}
$M^{p,q}_T$ is a Banach space with norm 
 \begin{align*}
\|(u,Q,P)\|_{M^{p,q}_T}&:=\|u\|_{L^\infty([0,T],D^{1-\frac{1}{p},p}_{A_q})}+\|u\|_{L^p(0,T;W^{2,q}(U)\cap W^{1,q}_0(U))}+\|\partial_tu\|_{L^p(0,T;L^q(U))}\\
&\ \ +\|Q\|_{L^\infty([0,T];B^{3-\frac{2}{p}}_{q,p}))}+\|Q\|_{L^p(0,T;W^{3,q}(U))}+\|\partial_tQ\|_{L^p(0,T;W^{1,q}(U))}\\
&\ \ +\|P\|_{L^p(0,T;W^{1,q}(U))}.
\end{align*}
\end{definition}
In the above definition, the space $D^{1-\frac{1}{p},p}_{A_q}$ is some fractional domain of the  Stokes operator in $L^q$ , for details see  \cite{MR2258416}, Section 2.3.
In short, this space consists of vectors whose $2-\frac{2}{p}$ derivatives are in $L^q$, divergence free, and vanish on the boundary of the domain.
Moreover, we have the following embedding relation
\begin{equation*}
  D^{1-\frac{1}{p},p}_{A_q}\hookrightarrow B^{2(1-\frac{1}{p})}_{q,p}\cap L^q.
\end{equation*}
Besov spaces on  bounded domains $U$ can be defined as interpolation spaces as in \cite{MR0482275},
\begin{equation*}
  B^{2(1-\frac{1}{p})}_{q,p}=(L^q,W^{2,q})_{1-\frac{1}{p},p}.
\end{equation*}
Notice that for the incompressible system, the pressure is still a solution up to a difference by a constant, so the condition on  $ P$ in Definition
\ref{CH1} is satisfied automatically, if we replace the pressure by $P-\frac{1}{|U|}\int_UPdx.$
\begin{definition}\label{DE1}
Let $T\in(0,+\infty]$. The triplet $(u,Q,P)$ is a strong solution to \eqref{cz22}, \eqref{cz33} on $U\times[0,T)$ if $(u,Q,P)$  belongs to $M^{p,q}_t$ whenever $t<T$, it satisfies the equation \eqref{cz22} in distribution sense on $U\times(0,T)$ and the initial boundary conditions \eqref{cz33}.
\end{definition}
 Recently, Xiao \cite{MR3582195}   proved the existence of a global strong solution for  small initial data
 if one  of the following conditions\par \noindent
$ {\bf (A_1)}:$\,\,\,\,\,\,\,\,\,\,\,\,\,\,\,\,\,\,\,\,\,\,\,  $  ac>\frac{9}{16}b^2,\,  a>0, \, c>0,$ \,\,\,\,\,\,\,\,\,\,\,\,\,\,\, or \hfill
${\bf (A_2)}:$\,\,\,\,\,\,\,\,\,\,\,\,\,\,\,\,\,\,\,\,\,\,\,\,\,\,\,\,$ a=0$,\,\,\,\,\,\,\,\,\,\,\,\,\,\,\,\,\,\,\,\par\noindent
is satisfied. Specifically, Xiao's result is the following;
 \par
\begin{theorem}\label{xiao} \cite{MR3582195}
 Let $U$ be a bounded domain in $\R^3$ with regular boundary. Assume that $1<p<\infty, q>3, u_0\in D^{{1-\frac{1}{p}},p}_{A_q}, Q_0\in {{{B}^{3-\frac{2}{p}}_{q,p}}}\cap W^{1,q}.$ Then
\par\noindent
1. there exists $0<\epsilon_0\ll 1$ such that if
\begin{equation*}
  \|Q_0\|_{L^\infty}\leq \epsilon_0,
\end{equation*}
then there exists some $T>0$ such that the  system \eqref{cz22} and \eqref{cz33}
 has a unique local strong
solution $(u,Q,P)\in M^{p,q}_t$ in $U\times [0,T)$ for any $0<t<T.$
\par\noindent
2. If
$({\bf A_1})$  or $({\bf A_2})$
is satisfied, there exists a positive $\delta_0,$ such that if the initial data satisfies
\begin{equation*}
  \|u_0\|_ {D^{{1-\frac{1}{p}},p}_{A_q}}\leq \delta_0,\,\, \|Q_0\|_{{{{B}^{3-\frac{2}{p}}_{q,p}}}\cap W^{1,q}}\leq \delta_0,
\end{equation*}
then \eqref{cz22} admits a unique global strong solution $(u,Q,P)\in M^{p,q}_{T}$ in $U\times [0,T)$ for all $T>0.$
\end{theorem}
We point out  that $W^{1,q}\hookrightarrow L^\infty,$ if $q>3.$ i.e.
there exists a positive constant $C>0$ such that  $\|Q_0\|_{L^\infty}\leq C\|Q_0\|_{W^{1,q}}.$
This inequality implies $\|Q_0\|_{L^\infty}\leq \epsilon_0$ if $\delta_0\leq C^{-1}\epsilon_0.$

In his paper,  inspired by the work of Danchin in \cite{MR2258416} and Hu and Wang in \cite{MR2628824},
   Xiao makes use of the maximal regularity of the Stokes operator and of the parabolic operator to get the existence of the local solution $(u,Q,P)$ by an approximation scheme.  Then to achieve the existence of the global strong solution,  he uses a bootstrap argument to prove
   a uniform boundedness of the norm in the space $M^{p,q}_T$ of the local solution $(u,Q,P).$
   The key tool of the proof is an exponential decay of the norm $\|Q\|_{L^q},$ with $2\leq q\leq+\infty$
   which is valid under the assumption $({\bf A_1}).$

Our first result is an improvement of Xiao's global existence Theorem. Indeed, we prove that  the system \eqref{cz22} and \eqref{cz33}  has  a global  strong  solution   under the sole hypothesis $a\geq0$, removing   the more restrictive assumptions $({\bf A_1})$ and $({\bf A_2})$. As in Xiao's work,  the proof strategy involves showing  that the norm  $\|(u,Q,P)\|_{M_T^{p,q}}$ remains  uniformely bounded on any interval $[0,T]$.  The novelty is that we obtain the estimate thanks to a maximum principle (Lemma \ref{le10}), inspired  by a result in [7].

\par
Notice that $F(Q)=\frac{a}{2}|Q|^2+\frac{b}{3}trace[Q^3]+\frac{c}{4}|Q|^4, a>0$ has the following properties.
There exists a sufficiently small positive $r$ that depends on $(a,b,c)$ such that if
\begin{equation}\label{ce1}
  \mathcal{O}:=\{Q\in\R^{3\times3}_{sym,0}|~~ |Q|<r\},
\end{equation}
then $Q=0$ is the (unique) global minimum of $F$ in $\mathcal{O},$ i.e.
\begin{equation}\label{CWW7}
  F(0)=0,~~~F(Q)>0~~\text{for~~any~}~Q\in\mathcal{O}\setminus\{0\},
\end{equation}
and
\begin{equation}\label{mp6}
  \partial F(Q):Q\geq 0~~\text{for\,\, any\,}~Q\in\mathcal{O},
\end{equation}
\begin{equation}\label{mp6b}
F(Q)\geq \lambda |Q|^2, ~~\text{for\,\, some \,positive }~\lambda\,\,\,\text{and \,\, for\,\, any\,}~Q\in\mathcal{O}.
\end{equation}
Moreover, $F$ also satisfies
\begin{equation}\label{E2}
  c_1|Q|^2\leq\partial F(Q):Q\leq C|Q|^2, ~~\text{for~any}~~Q\in\mathcal{O}\setminus\{0\},
\end{equation}
for some $c_1\in (0,a)$ and some positive constant $C$.
\par
Using these properties we are able to prove the following Theorem.

 \noindent  {\sc { Theorem a}}.
 {\it Let $U$ be a bounded domain in $\R^3$ with regular boundary. Assume that $1<p<\infty,\ q>3,\ u_0\in D^{{1-\frac{1}{p}},p}_{A_q},\ Q_0\in {{{B}^{3-\frac{2}{p}}_{q,p}}}\cap W^{1,q},$ and $a\geq0$. Then
 there exists a positive $\delta_0,$ such that if the initial data satisfy
\begin{equation*}
  \|u_0\|_ {D^{{1-\frac{1}{p}},p}_{A_q}}\leq \delta_0, \,\,\,\|Q_0\|_{{{{B}^{3-\frac{2}{p}}_{q,p}}}\cap W^{1,q}}\leq \delta_0,
\end{equation*}
then the system \eqref{cz22} with  \eqref{cz33} admits  a unique global strong solution $(u,Q,P)\in M^{p,q}_{T}$ in $U\times [0,T)$ for all $T>0.$
Moreover, the solution satisfies
\begin{equation}\label{RC3}
  \|u(t)\|_{L^2}+\|Q(t)\|_{H^{1}}\leq C\ { exp}(-dt)~~\text{for any}~~ t\geq0~~~ \text{and some}~~~ d>0,
\end{equation}
and
\begin{equation}\label{RC4}
  \|Q(t)\|_{L^{q}}\leq \|Q_0\|_{L^q}\ { exp}(-c_1t),~~4\leq q\leq\infty,   \text{ and for any } t\geq 0.
\end{equation}
with $ c_1$  as  in  \eqref{E2}.}

\begin{remark}\label{r1}
Motivated by \cite{MR3486134}, we prove that if $a>0$ and the domain $U$ is bounded with regular boundary,
$\|Q\|_{L^q}, 1\leq q\leq+\infty$ has an exponential decay in time for
initial data sufficiently small.
Then by a bootstrap argument in a similar way as in \cite{MR3582195}, we can prove that the
norm $\|(u,Q,P)\|_{M^{p,q}_T}$ of the strong solution is uniformly bounded on any interval $[0,T]$
and so the solution   exists on any interval $[0,T], $ $\forall T>0.$ Moreover, we prove
also that the norms $\|u(t)\|_{L^2}$ and $\|Q(t)\|_{H^1}$
have an exponential decay in time.
\end{remark}
\par
Now, we give a result about continuous dependence of solutions.

\noindent {\bf {\sc Theorem b.}}
{\it
 Let $U$ be a bounded domain in $\R^3$ with regular boundary. Assume that $2\leq p<\infty, q>3, u_0\in D^{{1-\frac{1}{p}},p}_{A_q}, Q_0\in {{{B}^{3-\frac{2}{p}}_{q,p}}}\cap W^{1,q},$ and $\|Q_0\|_{W^{1,q}}\leq\sigma$ for some positive $\sigma,$ such that $\sigma<\epsilon_0,$ with  $\epsilon_0$ as in Theorem \ref{xiao}.
 If $(u^n_0,Q^n_0)_{n\in \mathbb{N}}$ tends to $(u_0,Q_0)$ in $D^{{1-\frac{1}{p}},p}_{A_q}\times {{{B}^{3-\frac{2}{p}}_{q,p}}}\cap W^{1,q},$
 then there exists a positive $T$ independent of $n$ such that the sequence of solutions $(u^n,Q^n)_{n\in \mathbb{N}}\in M^{p,q}_{T}$ tends to the solution $(u,Q)$ in $L^2$ norm. Moreover,
 \begin{equation*}
   \|Q^n-Q\|^2_{L^2}+\|\Delta (Q^n-Q)\|^2_{L^2}+\|u^n-u\|^2_{L^2}\rightarrow 0,~~\text{uniformly~~ for~~} t\in [0,T]~~\text{when}~~~ n\rightarrow\infty.
 \end{equation*}}

In the last part of  the paper, we consider the system \eqref{cz22}, \eqref{cz33} in the whole space $U=\R^3$ and we
study the partial regularity properties of  a  class of  weak solutions.
Results on the partial regularity of  the Q tensor model
 have already been proved by  Du, Hu and Wang \cite{MR4134151}. Indeed, they extend  the CKN Theorem \cite{MR673830, MR1488514} to the Q tensor system.
Later, Liu give some new partial regularity criterion in terms of the velocity gradient only in \cite{MR4305937}.
\par
 In this paper, we intend to give a new local regularity criterion assuming that $u(t)$ is small in $\dot{B}^{-1}_{\infty,\infty}$
 on a suitable interval  of time and the $L^\infty$ norm of $Q_0$ is small.  This result is inspired by a similar criterion for MHD systems in
 \cite{MR4310804} and to our knowledge,  it was not previously  known for  system (1.6). We remark that a key step in our proof will be the fact that the $L^\infty$ norm of $Q(t)$ remains small on some interval $[0,T]$, with $T>0$, if $\|Q_0\|_{L^\infty}$ is small enough, as it is proved in \cite{MR2864407}.
 \par
 In order to introduce our results, we will give the definition of a {\it suitable weak solution} (see also \cite{MR4134151} for a similar definition), which is a weak solution satisfying a local energy inequality (the definition of weak solution can be found  in \cite{MR4134151} , Definition 1.1).
 \begin{definition}\label{CH2}
A weak solution $(u,P,Q) $ with $u\in(L^\infty(0,\infty;L^2(\R^3))\cap L^2(0,\infty;H^1(\R^3))$, $P\in L^{\frac{3}{2}}(0,\infty;L^\frac{3}{2}(\R^3))$ and
$Q\in(L^\infty(0,\infty;H^1(\R^3))\cap L^2(0,\infty;H^2(\R^3)))$ of \eqref{cz22} is a suitable weak solution of \eqref{cz22}, if, in addition,
$(u,P,Q)$ satisfies the local energy inequality
\begin{equation}\label{SU}
\aligned
  \int_{\R^3}(|u|^2&+|\nabla Q|^2)\phi(x,t)dx+2\int_{\R^3\times[0,t]}(|\nabla u|^2+|\nabla^2 Q|^2)\phi(x,s)dxds\\
  &\leq\int_{\R^3\times[0,t]}(|u|^2+|\nabla Q|^2)(\partial_t\phi
  +\Delta\phi)(x,s)dxds\\
  &~~+\int_{\R^3\times[0,t]}[(|u|^2+2P)u\cdot\nabla\phi+2\nabla Q\otimes
  \nabla Q:u\otimes\nabla\phi](x,s)dxds\\&
  +2\int_{\R^3\times[0,t]}(\nabla Q\otimes\nabla Q-|\nabla Q|^2I_3):\nabla^2\phi
  (x,s)dxds\\
  &-2\int_{\R^3\times[0,t]}[Q,\Delta Q]\cdot u\otimes\nabla\phi(x,s)d\phi dxds\\
  &-2\int_{\R^3\times[0,t]}([w,Q]\cdot(\nabla Q\nabla\phi)+\nabla(\mathcal{L}[\partial F(Q)])\cdot\nabla Q\phi)(x,s)dxds
\endaligned
\end{equation}
 for any $0\leq\phi\in C^\infty_0(\R^3\times(0,t])$ and for any $t>0$.
 \end{definition}
 Our result is based on the following  regularity criterion established in \cite{MR4134151}:

\begin{theorem}
\label{Th1000}
  Let us assume that  $F(Q)=\frac{a}{2}|Q|^2-\frac{b}{3}trace[Q^3]+\frac{c}{4}|Q|^4, c>0$, $U=\R^3.$ There
 exists $\epsilon_1>0$ such that if $(u,Q): \R^3\times(0,\infty)\rightarrow\R^3\times S^{(3)}_0$ is a suitable
 weak solution of \eqref{cz22}, which satisfies, for $z_0:=(x_0,t_0)\in\Omega\times(0,\infty)$,
 \begin{equation}\label{du2}
   \limsup_{r\rightarrow0}\frac{1}{r}\int_{Q_r(z_0)}(|\nabla u|^2+|\nabla^2Q|^2)dxdt<\epsilon^2_1,
 \end{equation}
 where $Q_r(z_0):=B_r(x_0)\times[t_0-r^2,t_0]$,
 then $(u,Q)$ is smooth near $z_0.$
\end{theorem}

 \par\noindent
 \par
 We recall that in the Navier-Stoke framework, many regularity criterion have been proved in the so called "critical spaces" that are some
scaling invariant spaces, such as $ L^3(\R^3)$ and $\dot{B}^{-1}_{\infty,\infty}(\R^3)$ and they are related by the following embeddings:
 \begin{equation*}
   L^3(\R^3)\hookrightarrow L^{3,q}(\R^3)\hookrightarrow\dot{B}^{-1+\frac{3}{p}}_{p,q}(\R^3)\hookrightarrow \mathrm{BMO}^{-1}(\R^3)
   \hookrightarrow\dot{B}^{-1}_{\infty,\infty}(\R^3),~~(3<p,q<\infty).
 \end{equation*}

\noindent {\bf {\sc Theorem c.}}
{\it
Let us assume that   $u_0\in \mathbf{H},$
where $\mathbf{H}=\text{Closure~of}~ \{u\in C^\infty_0(\R^3,\R^3):div u=0\}~ \text{in} ~L^2(\R^3),$
and $Q_0\in H^1(\R^3)\cap L^\infty(\R^3).$  Let the triplet $(u, Q, P)$ be a suitable weak solution to the system \eqref{cz22} in $\R^3\times(0,\infty)$, with $c>0$, with the initial data $(u_0,Q_0)$.
Then, there exist two  positive constants $\varepsilon_0$ and $\varepsilon_1$ such that if $u
\in L^\infty(t_0-\tilde{\varrho}^2,t_0;\dot{B}^{-1}_{\infty,\infty}) $ satisfies, 
 for some $\tilde{\varrho}>0$,
\begin{equation}\label{cww1}
   \|u\|_{L^\infty(t_0-\tilde{\varrho}^2,t_0;\dot{B}^{-1}_{\infty,\infty})}<\varepsilon_0,
\end{equation}
and
\begin{equation}\label{YYY2}
 \|Q_0\|_{L^\infty(\R^3)}\leq \varepsilon_1.
\end{equation}
one gets that $(x_0,t_0)\in(\R^3\times\R^+)$ is a regular point (namely, the solution
  $(u,Q)$ is smooth
  near $(x_0,t_0)$).
}

\begin{remark}\label{r28}
Thanks to the Theorem 1.1 in \cite{MR4134151}, when the initial data $u_0\in \mathbf{H}$ and $Q_0\in H^1(\R^3, S^{(3)}_0)\cap L^\infty(\R^3,S^{(3)}_0),$ then there exists a global suitable weak solution $(u, Q, P).$
\end{remark}

\par
The bootstrap argument will be employed in our proof. A rigorous statement of the abstract bootstrap principle can be found in T. Tao's book (see \cite{TT2}).
\par
The paper is organized as follows. In section 2, we will give  some key Theorems about the
 maximal regularities of Stokes operator and parabolic operators, some Lemma  and maximum principle.
In section 3, based on this, we give the proof of the {\sc Theorem a, b}, and  {\sc c}.

\section{Preliminaries}
\par
In this section, we will recall some Lemmas and Theorems in \cite{MR1345385,MR2258416,MR1138838} which will  be used  in the following.
Define $\mathcal{W}(0,T):=W^{1,p}(0,T;W^{1,q})\cap L^p(0,T;W^{3,q}).$
Now, we recall the maximal regularity property  of the  parabolic and Stokes operators.
\begin{theorem}\cite{MR1345385}\label{Th3}
 Given $w_0\in B^{3-\frac{2}{p}}_{q,p}$ and $f\in L^p(0,T;W^{1,q}),$ the Cauchy problem
 \begin{equation*}
   \partial_tw-\Delta w=f,~~t\in(0,T),~~w(0)=w_0,
 \end{equation*}
 has a unique solution $w\in\mathcal{W}(0,T),$ and
 \begin{equation*}
   \|w\|_{\mathcal{W}(0,T)}\leq C(\|f\|_{L^p(0,T;W^{1,q})}+\|w_0\|_{B^{3-\frac{2}{p}}}),
 \end{equation*}
 where $C$ is independent of $w_0,$ $f$ and $T.$ Moreover, there exists a positive constant $c_0$ independent of
 $f$ and $T$ such that
 \begin{equation*}
   \|w\|_{\mathcal{W}(0,T)}\geq c_0\sup_{t\in(0,T)}\|w\|_{B^{3-\frac{2}{p}}_{q,p}}.
 \end{equation*}
 \end{theorem}
\begin{theorem}\label{Th4}
Let $U$ be a bounded domain with $C^{2+\epsilon}(\epsilon>0)$ boundary in $\R^3$ and $1<p,q<\infty.$
Assume that $u_0\in D^{{1-\frac{1}{p}},p}_{A_q}$ and $f\in L^p(0,\infty;L^q(U)).$ Then the system
\begin{equation*}
\begin{cases}
  \partial_tu-\Delta u+\nabla P=f,\\
  \int_{U}Pdx=0,~~\mathrm{div}\ u=0,\\
  u|_{\partial U}=0,~~~u|_{t=0}=u_0.
\end{cases}
\end{equation*}
has  a unique solution $(u, P)$ for any  $T>0$ which satisfies
\begin{equation*}
\aligned
  \|u(T)\|&_{D^{{1-\frac{1}{p}},p}_{A_q}}+\|(\partial_tu,\Delta u,\nabla P)\|_{L^p(0,T;L^q(U))}\\
  &\leq C(\|u_0\|_{D^{{1-\frac{1}{p}},p}_{A_q}}+\|f\|_{L^p(0,T;L^q(U))}),
\endaligned
\end{equation*}
where $C=C(p,q,U).$
\end{theorem}
\par
Next, we introduce some $L^\infty$ estimates for the derivatives of functions $f\in L^\infty(0,T; D^{{1-\frac{1}{p}},p}_{A_q})\cap L^p(0,T;W^{2,q})$.
 \begin{lemma}\cite{MR2258416}\label{le1}
 Let $1<p,q<\infty.$
 \par\noindent
(i) If $0<\frac{p}{2}-\frac{3p}{2q}<1,$ then it holds that
\begin{equation*}
  \|\nabla f\|_{L^p(0,T;L^\infty)}\leq CT^{\frac{1}{2}-\frac{3}{2q}}\|f\|^{1-\theta}_{L^\infty(0,T;D^{{1-\frac{1}{p}},p}_{A_q})}\|f\|^{\theta}_{L^p(0,T;W^{2,p})},
\end{equation*}
for some constant $C$ depending only on $U,p,q$ and $\frac{1-\theta}{p}=\frac{1}{2}-\frac{3}{2q}.$
\par\noindent
(ii) If $\frac{p}{2}-\frac{3p}{2q}>1,$ then it holds that
\begin{equation*}
  \|\nabla f\|_{L^p(0,T;L^\infty)}\leq CT^{\frac{1}{p}}\|f\|_{L^\infty(0,T;D^{{1-\frac{1}{p}},p}_{A_q})}
\end{equation*}
for some constant $C$ depending only on $U,p,q.$
\end{lemma}

\begin{lemma}\cite{MR2628824}\label{le2}
 Let $1<p,\ r<\infty.$ For $f\in W^{1,p}(0,T;L^r)$ and $f(0)\in L^r,$ then for any $t\in[0,T],$
 \begin{equation*}
   \|f\|_{L^\infty(0,t; {L^r})}\leq C(\|f(0)\|_{L^r}+\|f\|_{W^{1,p}(0,t; L^r)}),
 \end{equation*}
 for some constant $C$ independent of $f$ and $T.$
\end{lemma}
Before proving the existence of global strong solutions, we give a key lemma which plays an important role in completing the bootstrap argument.
Motivated by the results in \cite{MR3486134}, we will  prove that the local solution $(u,P,Q)\in M^{p,q}_T$ for some $T>0$
also satisfies the following  interesting maximum principle.

\begin{lemma}\label{le10}
 Assume that    $T\in(0, +\infty)$ and that  $(u,Q,P)\in M^{p,q}_T$,\ $1<p<\infty,\ q>3$  is a  local strong solution to  system \eqref{cz22} with  \eqref{cz33} and with $a>0$.
If the initial data  $Q_0$ satisfy  $\|Q_0\|_{L^\infty}<r$ with $r$ as in  \eqref{ce1}, then
 \begin{equation}\label{mp4}
  \|Q(t)\|_{L^\infty}\leq \|Q_0\|_{L^\infty},
 \end{equation}
for any $t\in(0,T)$.
Moreover, if  $\|Q_0\|_{L^\infty}\leq\frac{r}{2}$, we have the better estimate
\begin{equation}\label{mp5}
  \|Q(t)\|_{L^{q}}\leq exp(-c_1t)\|Q_0\|_{L^{q}},~~~4\leq q\leq\infty,\  \ \text{for any }t\in (0,T)
\end{equation}
where  $c_1$ is the positive constant defined in \eqref{E2}.
\end{lemma}
\begin{proof}
{To prove the inequality \eqref{mp4}, we use the same idea as in \cite{MR3486134}.}
Firstly, we rewrite the second equation of the  system \eqref{cz22} as
\begin{equation}\label{111}
  \partial_tQ+u\cdot\nabla Q-\Delta Q=-\mathcal{L}[\partial F(Q)]+w(u)Q-Qw(u),~~Q(0,\cdot)=Q_0,
\end{equation}
where $\mathcal{L}[\partial F(Q)]=aQ-b\left[Q^2-\frac 13 Id \ trace (Q^2)\right]+cQ trace(Q^2)$.
Note that $Q\in{C([0,T];B^{3-\frac{2}{p}}_{q,p})}$ and so by
the following embeddings (see \cite{MR0482275})
\begin{equation*}
  \|Q\|_{L^\infty}\leq\|Q\|_{B^{0}_{\infty,1}}\leq\|Q\|_{B^{3-\frac{2}{p}-\frac{3}{q}}_{\infty,p}}\leq\|Q\|_{B^{3-\frac{2}{p}}_{q,p}},
\end{equation*}
 the local solution $Q\in C([0,T];L^\infty)$.
So we can
take the scalar product of \eqref{111} with $2G'(|Q|^2)Q,$ where $G\in C^2([0,\infty))$ is a suitable chosen function, and integrate over the domain $U$ to obtain with some computations:
\begin{equation}\label{YJ2}
  \frac{d}{dt}\int_{U}G(|Q|^2)dx+\int_{U}[2G'(|Q|^2)|\nabla Q|^2+G''(|Q|^2)|\nabla|Q|^2|^2]dx=-2\int_{U}G'(|Q|^2)\partial F(Q): Qdx.
\end{equation}
Indeed,  the equality  \eqref{YJ2} is a consequence of
 the property $\mathrm{div}~u=0$ and
 the {following equality}
\begin{equation}\label{YJ100}
\aligned
  (w(u)Q-Qw(u))&:Q=-(w(u)Q+Qw(u)):Q+2w(u) Q:Q=-(w(u)Q+Qw(u))_{\alpha\beta}Q_{\alpha\beta}\\
  &+2w(u)_{\alpha z}Q_{z\beta}Q_{\alpha\beta}=0,
\endaligned
\end{equation}
which  is a consequence of the fact that the terms $(w(u)Q+Qw(u))_{\alpha\beta}Q_{\alpha\beta}$ and  $(w(u)Q+Qw(u))_{\beta\alpha}Q_{\beta\alpha}$
 appear in pairs and  $w(u)Q+Qw(u)$ is an antisymmetric  matrix while $Q$ is a symmetric  matrix.
 Similarly, the term $w(u)_{\alpha z}Q_{z\beta}Q_{\alpha\beta}$ and $w(u)_{z\alpha }Q_{\alpha\beta}Q_{z\beta }$ appear in pairs and $w(u)$ is an antisymmetric matrix.
In order to get  \eqref{YJ2}, we have  also used the property
$$
\int_{U}G'(|Q|^2){\mathcal L}\partial F(Q): Qdx=\int_{U}G'(|Q|^2)\partial F(Q): Qdx,
$$
since
 \begin{equation*}
   2\int_{U}G'(|Q|^2)Q_{\alpha\beta}\frac{1}{3}trace(\partial F(Q))I_{\alpha\beta}dx=2\int_{U}G'(|Q|^2)Q_{\alpha\alpha}\frac{1}{3}trace(\partial F(Q))dx=0.
 \end{equation*}

Now take $G$ in \eqref{YJ2} such that
\begin{equation}\label{E1}
  G(z)=0~~\text{for}~z\in[0,\|Q_0\|^2_{L^\infty}],~~G'\geq 0,~~G''\geq 0,~~G(z)>0~~~~\text{for}~z>\|Q_0\|^2_{L^\infty}
\end{equation}
{Since $\|Q_0\|_{L^\infty}<r$ and  $\mathcal{O}$ is open, there exists  $T_1>0$ small enough  such that $Q(t,\cdot)\in\mathcal{O}$ for any $0\leq t\leq T_1$.}
 Since $G'\geq0,$ $G''\geq0$  and using \eqref{mp6}, then \eqref{YJ2} implies
 \begin{equation*}
   \int_{U}G(|Q|^2)(t)dx\leq\int_{U}G(|Q_0|^2)dx=0~~\text{for~~all~}0\leq t\leq T_1.
 \end{equation*}
 So $G(|Q|^2(t,x))=0$ for a.e. $x\in U$ and  thanks to \eqref{E1}, we get
\begin{equation}\label{YJ1}
  \|Q(t)\|_{L^\infty}\leq \|Q_0\|_{L^\infty}~~\text{for~~all~}0\leq t\leq T_1.
\end{equation}
Now, we let us define  $\tilde{T}:=\sup\{t>0, \|Q(t)\|_{L^\infty}\leq \|Q_0\|_{L^\infty}\}.$
Then we can prove that $\tilde{T}=T$ by a contradiction argument. If $\tilde{T}<T,$ note $\|Q(\tilde{T})\|_{L^\infty}\leq\|Q_0\|_{L^\infty}<r,$
and using the same argument as above, we  get
 \begin{equation*}
   \|Q(t)\|_{L^\infty}\leq\|Q_0\|_{L^\infty},~\text{for}~{\tilde{T}\leq t\leq \tilde{T}+\delta},
 \end{equation*}
 for some positive $\delta.$ But this  contradicts the definition of the maximal time span $\tilde T$. Therefore, we get  $\tilde{T}=T$
and
\begin{equation}\label{YJ1}
  \|Q(t)\|_{L^\infty}\leq \|Q_0\|_{L^\infty}~~\text{for~~all~} 0\leq t\leq T.
\end{equation}
We could  prove the inequality \eqref{mp5} using the inequality \eqref{mp4}.
 But we prefer to prove it indipendently by a bootstrap argument.
Firstly, we prove
\begin{equation}\label{RC1}
  \|Q(t)\|_{L^\infty}\leq \frac{r}{2},~~~\forall t\in[0,T].
\end{equation}
by bootstrap argument. We assume that
$\|Q(t)\|_{L^\infty}\leq r,~~~\forall t\in[0,T]$ and we prove that this implies $\|Q(t)\|_{L^\infty}\leq \frac{r}2,~~~\forall t\in[0,T]$, if $\|Q_0\|_{L^\infty}\leq \frac{r}2$.
Let $G(z)=z^q, q\geq2$ in \eqref{YJ2}. By  simple computations, we get
\begin{equation}\label{Eq}
  \frac{d}{dt}\int_{U}|Q|^{2q}dx+\int_{U}[2q|Q|^{2{(q-1)}}|\nabla Q|^2+q(q-1)|Q|^{2(q-2)}|\nabla|Q|^2|^2]dx=-2\int_{U}q|Q|^{2(q-1)}\partial F(Q):Qdx
\end{equation}
Note that the term $\int_{U}[2q|Q|^{2{(q-1)}}|\nabla Q|^2+q(q-1)|Q|^{2(q-2)}|\nabla|Q|^2|^2]dx$ is positive  and since  $\|Q(t)\|_{L^\infty}\leq \frac{r}2$  then $Q(t)\subset\mathcal{O}$, $\forall t\in[0,T]$, and so thanks to \eqref{E2} we have $$\int_{U}q|Q|^{2(q-1)}\partial F(Q):Qdx\geq c_1\int_{U}q|Q|^{2q}dx.$$
 So we get from \eqref{Eq}
\begin{equation*}
  \frac{d}{dt}\int_{U}|Q|^{2q}dx\leq -c_1\int_{U}q|Q|^{2q}dx.
\end{equation*}
It follows by Gronwall's inequality that
\begin{equation}\label{mp1}
  \|Q(t)\|_{L^{2q}}\leq exp(-c_1t)\|Q_0\|_{L^{2q}},~~~~q\geq2.
\end{equation}
{Notice the constant $c_1$ is independent of $q,$ and so when we let   $q\rightarrow\infty,$ we have
\begin{equation}\label{mp2}
  \|Q(t)\|_{L^{\infty}}\leq exp(-c_1t)\|Q_0\|_{L^{\infty}},
\end{equation}}
 for any $t\in (0,T)$.
Therefore if  we choose $\|Q_0\|_{L^{\infty}}\leq\frac{r}{2}$, we get $\|Q(t)\|_{L^{\infty}}\leq\frac{r}{2},$ for any $t\in (0,T)$.

By bootstrap this proves the claim \eqref{RC1}. Moreover, if $\|Q_0\|_{L^{\infty}}\leq\frac{r}{2}$ also the inequalities \eqref{mp1} and \eqref{mp2} hold and so the inequality \eqref{mp5} is proved.
So  the proof of the Lemma is completed.
\end{proof}
\par
\begin{remark}\label{r10}
We remark that  in the proof of  the Lemma \ref{le10}, we prove the inequality \eqref{mp4} only using the
properties \eqref{CWW7} and \eqref{mp6}.
\end{remark}
\section{Global existence of  a strong solution}
\par\noindent
{\textbf{Proof of {\sc Theorem a.}}
\par
1. {\it  Local existence of strong solution}.
\par\noindent
By  Theorem \ref{xiao}  if $\|Q_0\|_{L^\infty}\leq \epsilon_0$ there exist some positive time $T$ and a
 local solution $(u,Q,P)\in M^{p,q}_{t}$  for any $0<t<T$.
\par\noindent
2. Decay of the solution and  energy inequalities.
\par\noindent
In order to establish the existence of   the global solution, we employ a bootstrap argument. Let $(u,Q,P)$ be the local solution that belongs to  $M^{p,q}_{t}$ for any $0<t<T$ and for some positive $T$. As in the paper \cite{MR3582195},
we define
\begin{equation}\label{ene}
\begin{aligned}
 { H}(t):=&\|u\|_{L^\infty(0,t;D^{{1-\frac{1}{p}},p}_{A_q})}+\|u\|_{L^p(0,t;W^{2,q})}+\|\nabla P\|_{L^p(0,t,L^q)}\\
  &+\|\partial_tu\|_{L^p(0,t;L^q)}+\|Q\|_{L^\infty(0,t;B^{3-\frac{2}{p}}_{q,p})}+\|Q\|_{\mathcal{W}(0,t)},
\end{aligned}
\end{equation}
and
\begin{equation*}
{H}_0:=\|u_0\|_{D^{{1-\frac{1}{p}},p}_{A_q}}+\|Q_0\|_{B^{3-\frac{2}{p}}_{q,p}\cap W^{1,q}},
\end{equation*}
where we recall that
\begin{equation*}
  \mathcal{W}(0,t)=W^{1,p}(0,t;W^{1,q})\cap L^p(0,t;W^{3,q}).
\end{equation*}
If we are able to prove that $H(t)$ is uniformly bounded for any $t\in (0,T)$ then the
 solution $(u,Q,P)$ can be extended to any interval $ (0,T)$. Now, we  prove that ${H}(t)$ is uniformly  bounded on $(0,T)$.
By  Theorems \ref{Th3} and \ref{Th4} and thanks to the Poincar\'e inequality, we have
\begin{equation*}
\aligned
  {H}(t)\leq &C\Big({H}_0+\|u\cdot\nabla u\|_{L^p(0,t;L^q)}+\|\mathrm{div}(\nabla Q\otimes\nabla Q)\|_{L^p(0,t;L^q)}\\
  &+\|\mathrm{div}(\Delta QQ-Q\Delta Q)\|_{L^p(0,t;L^q)}+\|u\cdot\nabla Q\|_{L^p(0,t;W^{1,q})}\\
  &+\|Q\omega-\omega Q\|_{L^p(0,t;W^{1,q})}+\|-aQ\|_{L^p(0,t;W^{1,q})}\\
  &+\|b[Q^2-\frac{Id}{3}trace(Q^2)]-cQtrace(Q^2)\|_{L^p(0,t;W^{1,q})}\Big)\\
  &:=C\big({H}_0+K_1+K_2+K_3+K_4+K_5+K_6+K_7\big).
\endaligned
\end{equation*}
As in the paper \cite{MR3582195}, we can estimate $K_1-K_5$ and $K_7$ in the following way:
\begin{equation*}
\aligned
  K_1&\leq \|u\|_{L^\infty(0,t;L^q)}\|\nabla u\|_{L^p(0,t;L^\infty)}\\
  &\leq C\|u\|_{L^\infty(0,t;D^{1-\frac{1}{p},p}_{A_q})}\|u\|_{L^p(0,t;W^{2,q})}\\
  &\leq C{H}^2(t).
\endaligned
\end{equation*}
Moreover,
\begin{equation*}
\aligned
  K_2&\leq\|\nabla Q\|_{L^\infty(0,t;L^q)}\|\nabla^2Q\|_{L^p(0,t;L^\infty)}\\
  &\leq\| Q\|_{L^\infty(0,t;W^{1,q})}\|Q\|_{L^p(0,t;W^{3,q})}\\
  &\leq C(\|Q_0\|_{W^{1,q}}+\|Q\|_{W^{1,p}(0,t;W^{1,q})})\|Q\|_{\mathcal{W}(0,t)}\\
  &\leq C({H}_0+{H}(t)){H}(t),
\endaligned
\end{equation*}
where in the third inequality we used the Lemma \ref{le2}.
Then,  denoting by $\nabla^2 Q$ and  $\nabla^3 Q$, respectively, the set of all the second and third order  partial derivatives, by Holder inequality and Lemma  \ref{le2}, we get
\begin{equation*}
\aligned
K_3
&\leq C\|Q\|_{L^\infty(0,t;L^\infty)}\|\nabla^3Q\|_{L^p(0,t;L^q)}+C\|\nabla Q\|_{L^\infty(0,t;L^q)}\|\nabla^2Q\|_{L^p(0,t;L^\infty)}\\
&\leq C(\|Q_0\|_{W^{1,q}}+\|Q\|_{W^{1,p}(0,t;W^{1,q})})\|Q\|_{L^p(0,t;W^{3,q})}\\
&\leq C({H}_0+{H}(t)){H}(t);
\endaligned
\end{equation*}
\begin{equation*}
\aligned
K_4&\leq \|u\|_{L^\infty(0,t;L^q)}\|\nabla Q\|_{L^p(0,t;W^{1,\infty})}+\|\nabla u\|_{L^p(0,t;L^\infty)}\|\nabla Q\|_{L^\infty(0,t; L^q)}\\
&\leq C\|u\|_{L^\infty(0,t;D^{{1-\frac{1}{p}},p}_{A_q})}\|Q\|_{L^p(0,t;W^{3,q})}\\
&~~~~~+C\|u\|_{L^p(0,t;W^{2,q})}(\|Q_0\|_{W^{1,q}}+\|Q\|_{W^{1,p}(0,t; W^{1,q})})\\
&\leq C({H}_0+{H}(t)){H}(t);
\endaligned
\end{equation*}
\begin{equation*}
\aligned
K_5&\leq C\|\nabla u\|_{L^p(0,t;L^\infty)}\|Q\|_{L^\infty(0,t;W^{1,q})}+\|\nabla^2u\|_{L^p(0,t; L^q)}\|Q\|_{L^\infty(0,t;L^\infty)}\\
&\leq C({H}_0+{H}(t)){H}(t);
\endaligned
\end{equation*}
\begin{equation*}
\aligned
K_7&\leq C\|Q\|_{L^\infty(0,t;L^\infty)}\|Q\|_{L^p(0,t;W^{1,q})}+\|Q\|^2_{L^\infty(0,t;L^\infty)}\|Q\|_{L^p(0,t;W^{1,q})}\\
&\leq C(H^2_0+H^3_0+H^2(t)+H^3(t)).
\endaligned
\end{equation*}
It remains to deal with the term $K_6$. We stress that under the
 assumption $a>0,$  we cannot adopt the same  method as in paper \cite{MR3582195}.
Let  us assume that $\delta_0$ in the hypothesis of Theorem \ref{xiao} is small enough such that
$\|Q_0\|_{L^\infty}\leq C\|Q_0\|_{W^{1,q}} \leq\frac{r}{2},$ with $r$ as in \eqref{ce1}.
So thanks to Lemma \ref{le10}  the solution $Q(t)$ satisfies the inequality \eqref{mp5} on $[0,T]$.
Since  $q>3$ and  the domain is bounded, thanks to the inequality \eqref{mp2}, the embeddings $W^{1,q}(U)\hookrightarrow L^\infty(U)\hookrightarrow L^q(U),\ q>3$ and the interpolation Gagliardo-Niremberg inequality we get
\begin{equation*}
\aligned
  K_6:=\|-aQ\|&_{L^p(0,t;W^{1,q})}\leq \|Q\|_{L^p(0,t;L^q)}+\|\nabla Q\|_{L^p_t(L^q)}\\
  &\leq C\big(\int^t_0\|Q(s,\cdot)\|^p_{L^\infty}ds\big)^\frac{1}{p}+
  C\|Q\|^\frac{1}{2}_{L^p(0,t;L^q)}\|\nabla^2Q\|^{\frac{1}{2}}_{L^p(0,t;L^q)}\\
  &\leq C\|Q_0\|_{L^\infty}\big(\int^t_0\exp(-pc_1s)ds\big)^\frac{1}{p}+C\big(\|Q_0\|_{L^\infty}(\int^t_0\exp(-pc_1s)ds)^\frac{1}{p}\big)^\frac{1}{2}H^\frac{1}{2}(t)\\
  &\leq CH_0+CH_0^\frac{1}{2}H^\frac{1}{2}(t)\\
  \endaligned
  \end{equation*}
  Now, by Young inequality, for any $\varepsilon>0$
 \begin{equation}\label{Ye}
 \aligned
  &CH_0+CH_0^\frac{1}{2}H^\frac{1}{2}(t)\\
  &\leq CH_0+C\epsilon^{-1}H_0+C\epsilon H(t).
\endaligned
\end{equation}
Choosing $\epsilon$ small enough in \eqref{Ye} and collecting  the estimates on  the terms $K_1, K_2, \cdots, K_7$, one obtains
\begin{equation}\label{lb99}
  {H}(t)\leq C({\hat{H}}_0+{H}^2(t)+{H}^3(t)),
\end{equation}
where ${\hat{H}}_0:={H}_0+{H}^2_0+{H}^3_0$ and  $C$ is a constant independent of ${\hat{H}}_0$ and of time $t$.
The bootstrap argument starts with the following assumption:
 ${H}(t)$ is bounded uniformly for $t\in[0,T]$,  i.e.
\begin{equation}\label{RC2}
  {H}(t)\leq 4C{\hat{H}}_0,~~\forall t\in[0,T],
\end{equation}
for some positive constant $C$.
If the initial data are  such that $C\|Q_0\|_{W^{1,q}} \leq\frac{r}{2}$ and
\begin{equation*}
  \max\{16C^2{\hat{H}}_0+64C^3\hat{H}^2_0\}\leq1
\end{equation*}
then it follows from the inequality \eqref{lb99} and \eqref{RC2} that
\begin{equation*}
  {H}(t)\leq 2C{\hat{H}}_0,~~\forall t\in[0,T].
\end{equation*}
So we complete the proof of the bootstrap assumption \eqref{RC2}. Moreover, we get that the
solution is global and the decay estimates  \eqref{RC4}.
\par\noindent
3. The decay of  $\|u\|_{L^2}$ and $\|Q\|_{H^1}$.
\par\noindent
Now, we will prove the inequality \eqref{RC3}.
As in \cite{MR3486134} choosing  $G(z)=\sqrt{z},$ by a direct computation one can check that the sum $2G'(|Q|^2)|\nabla Q|^2
+G''(|Q|^2)\big|\nabla|Q|^2\big|^2$
is nonnegative.
Combining the fact that  $Q(t,\cdot)\subset\mathcal{O}$  for any $t\geq0$ and using   \eqref{mp6}, we get $-2\int_{U}G'(|Q|^2)\partial F(Q): Q~dx\leq 0.$
From  \eqref{YJ2} it follows  that
\begin{equation}\label{YJ3}
  \|Q(t)\|_{L^1}\leq \|Q_0\|_{L^1},~~~\text{for~~all~}t>0.
\end{equation}
Finally, taking  $G(z)=z$ in \eqref{YJ2},  we obtain
\begin{equation*}
  \frac{d}{dt}\int_{U}\frac{1}{2}|Q|^2~dx+\int_{U}|\nabla Q|^2~dx\leq -\int_{U}\partial F(Q):Qdx~~\text{for~all~}t>0.
\end{equation*}
By \eqref{E2}, since $Q(t,\cdot)\subset\mathcal{O}$  for any $t\geq 0$, we have
\begin{equation*}
  \partial F(Q):Q\geq c_1|Q|^2.
\end{equation*}
Integrating in time, we have
\begin{equation*}
  \frac{1}{2}\|Q(t)\|^2_{L^2}+\int^t_0\|\nabla Q\|^2_{L^2}~ds\leq -c_1\int^t_0\|Q(t)\|^2_{L^2}ds+\frac{1}{2}\|Q_0\|^2_{L^2}
\end{equation*}
Using the Gronwall's inequality, we obtain
\begin{equation}\label{YJ10}
  \|Q(t)\|^2_{L^2}\leq \|Q_0\|^2_{L^2}exp(-2c_1t)~~\text{for~all}~t\geq0.
\end{equation}
\par\noindent
On the other hand, thanks to  the regularity of the local solution $(u,Q),$ using the method in \cite{MR2864407} and the boundary condition $(u,\partial_\nu Q)\mid_{\partial U}=(0,\mathcal{O}),$  we get the following a priori estimate
\begin{equation*}
\aligned
  \int_{U}\Big [\frac{1}{2}|u|^2&+\frac{1}{2}|\nabla Q|^2+F(Q)\Big ](t,x)dx+\int^t_s\int_{U}\left[|\nabla u|^2+|\Delta Q-\mathcal{L}[\partial F(Q)]|^2\right](\sigma,x)~ dx~ d\sigma\\
  &\leq\int_{U}\Big[\frac{1}{2}|u|^2+\frac{1}{2}|\nabla Q|^2+F(Q)\Big](s,x)~dx
\endaligned
\end{equation*}
for all $t>s$ and a.e. $s\in[0,\infty)$ including $s=0.$
We can rewrite the above inequality into the  differential version:
\begin{equation*}
  \frac{d}{dt}\int_{U}\left[\frac{1}{2}|u|^2+\frac{1}{2}|\nabla Q|^2+F(Q)\right]dx+\int_{U}\left[|\nabla u|^2+|\Delta Q-\mathcal{L}[\partial F(Q)]|^2\right](t,x)~dx\leq 0.
\end{equation*}
Using the properties \eqref{mp6} and \eqref{E2}, we prove that  the energy inequality reads as follows:
\begin{equation}\label{YJ5}
  \aligned
  \frac{d}{dt}\int_{U}\Big[\frac{1}{2}|u|^2&+\frac{1}{2}|\nabla Q|^2+F(Q)\Big](t,x)~dx\\
  &+C\int_{U}[\frac{1}{2}|\nabla u|^2+\frac{1}{2}|\nabla Q|^2+F(Q)](t,x)~dx\leq0,~~C>0.
\endaligned
\end{equation}
Indeed, this is a consequence of the following inequality
\begin{equation}\label{s3}
  \int_{U}|\Delta Q-\mathcal{L}[\partial F(Q)]|^2~dx\geq C\int_{U}[\frac{1}{2}|\nabla Q|^2+F(Q)]~dx,~~C>0.
\end{equation}
To prove \eqref{s3}, let
\begin{equation*}
  -\Delta Q+\mathcal{L}[\partial F(Q)]=Z,~~Z:=-u\cdot\nabla Q-Qw(u)+w(u)Q-\partial_tQ,
\end{equation*}
and observe that, thanks to \eqref{E2} and \eqref{mp6b}, we have
\begin{equation}\label{s1}
  \int_{U}\Big[\frac{1}{2}|\nabla Q|^2+F(Q)\Big]~dx\leq C\int_{U}[|\nabla Q|^2+\partial F(Q):Q]dx.
  \end{equation}
  Then,  by integration by parts and since $\partial F(Q):Q={\mathcal L}(\partial F(Q)):Q$, we get
  \begin{equation}\label{s2}
  \aligned
  \int_{U}[|\nabla Q|^2+\partial F(Q):Q]dx & =C\int_{U}Z:Qdx\\
  &\leq C(\varepsilon)\|Z\|^2_{L^2}+\varepsilon\|Q\|^2_{L^2},
\endaligned
\end{equation}
where in the last inequality we use  H\"older and Young  inequalities with $\varepsilon>0$.
Finally, collecting the inequalities \eqref{s1} and \eqref{s2}, we obtain
\begin{equation}
 \int_{U}\Big[\frac{1}{2}|\nabla Q|^2+F(Q)\Big]~dx\leq C(\varepsilon)\|Z\|^2_{L^2}+\varepsilon\|Q\|^2_{L^2}
\end{equation}
and since,   for $\varepsilon>0$ sufficiently small,
the term $\varepsilon\|Q\|^2_{L^2}$ can be "absorbed" by $\int_{U}F(Q)~dx$, thanks to \eqref{mp6b} and  the  inequality \eqref{s3} is proved.

Using  Poincare inequality in \eqref{YJ5}, we obtain
\begin{equation}\label{YJ7}
  \aligned
  \frac{d}{dt}\int_{U}[\frac{1}{2}|u|^2&+\frac{1}{2}|\nabla Q|^2+F(Q)](t,\cdot)dx\\
  &\leq -C\int_{U}[\frac{1}{2}|u|^2+\frac{1}{2}|\nabla Q|^2+F(Q)](s,\cdot)dx,~~C>0.
\endaligned
\end{equation}
It follows from  Gronwall's inequality that
\begin{equation*}
  \int_{U}\Big[\frac{1}{2}|u|^2+\frac{1}{2}|\nabla Q|^2+F(Q)\Big](t,x)dx\leq (\|u_0\|^2_{L^2}+\|\nabla Q\|^2_{L^2}+\|F(Q_0)\|^2_{L^2})exp(-Ct)
\end{equation*}
Then, using \eqref{mp6b},
one has
\begin{equation}\label{YJ11}
  \|u(t)\|_{L^2}+\|Q(t)\|_{H^1}\leq C(\|u_0\|_{L^2}+\|\nabla Q_0\|_{L^2}+\|F(Q_0)\|_{L^2})exp(-Ct)
\end{equation}
On the other hand, the inequality \eqref{RC4} has been proved in Lemma \ref{le10}.
Thus we completed the proof of  {\sc Theorem a}.
\par

\par\noindent
{\bf Proof of {\sc Theorem b}.}
\begin{proof}
Now, we prove the continuous dependence of the solutions with respect to the initial data.
Since $(u^n_0,Q^n_0)\rightarrow(u_0,Q_0)$ in $D^{{1-\frac{1}{p}},p}_{A_q}\times {{{B}^{3-\frac{2}{p}}_{q,p}}}\cap W^{1,q}$
and $\|Q_0\|_{W^{1,q}}<\sigma,$ then
$\|Q^n_0\|_{L^\infty}\leq\|Q^n_0\|_{W^{1,q}}\leq\sigma,$
for large $n.$ So, thanks to Theorem \ref{xiao}, there exists a sequence of local solutions $(u^n,Q^n)$
on $(0,T^n),$ for some positive $T^n.$ Moreover, since $\|Q^n_0\|_{L^\infty}\leq\sigma$ for $n$  large,
it is possible to choose a common existence time $T>0.$
As in \eqref{ene}, we introduce the energy of any solution $(u^n,Q^n,P^n)$
\begin{equation*}
 \aligned
  {H}_n(t):=&\|u^n\|_{L^\infty(0,t;D^{1-\frac{1}{p},p}_{A_q})}+\|\partial_tu^n \|_{L^p(0,t;L^q)}+\|u^n \|_{L^p(0,t;L^q)}\\ &+\|\nabla P^n \|_{L^p(0,t;L^q)} +\|Q^n\|_{L^\infty(0,t;B^{3-\frac{2}{p}}_{q,p})}
  +\|Q^n\|_{{\mathcal W}{(0,t)}},
\endaligned
\end{equation*}
and, by similar estimates as in the proof of {\sc Theorem a}, we get  that $H_n(t)$ and $H(t)$,
 are bounded on $[0,T]$ uniformly in $n.$
Next, we will show the solution $(u^n,Q^n)$ tends to the solution $(u,Q)$ in $L^2,$ uniformely for any $t\in [0,T]$.
Denote
\begin{equation*}
  \delta^n Q:=Q^n-Q,~~ \delta^n u:=u^n-u,~~ \delta^n P=P^n-P.
\end{equation*}
Then $(\delta^n Q,\delta^n u,\delta^n P)$ verifies the following system, at least in the weak sense, with
zero  boundary conditions
\begin{equation}\label{gs3}
\begin{cases}
\partial_t\delta^n Q+\delta^n u\cdot\nabla\delta^n Q-\delta^n w(u)\delta^n Q+\delta^n Q\delta^n w(u)\\
~~~=-\delta^n u\cdot\nabla Q-u\cdot\nabla \delta^n Q-Q\delta^n w(u)-\delta^n Qw(u)+\delta^n w(u)Q+w(u)\delta^n Q\\
~~~+\Delta\delta^n Q-a\delta^n Q+b[\delta^n QQ^n+\delta^nQ Q-\frac{Id}{3}tr(\delta^n QQ^n+\delta^n QQ)]\\
~~~-c[\delta^n Qtr(Q^n)^2+Qtr(\delta^n QQ^n+\delta^n QQ)],\\
\partial_t\delta^n u+\delta^n u\cdot\nabla\delta^n u-\Delta \delta^n u+\nabla \delta^n P\\
~~=-\delta^n u\cdot\nabla u-u\cdot\nabla\delta^n u-\mathrm{div}(\nabla\delta^n Q\otimes\nabla\delta^n Q)\\
~~~-\mathrm{div}(\nabla Q\otimes\nabla\delta^n Q+\nabla \delta^n Q\otimes\nabla Q)+\mathrm{div}(\delta^n Q\Delta\delta^n Q-\Delta\delta^n Q\delta^n Q)\\
~~~+\mathrm{div}(Q\Delta\delta^n Q+\delta^n Q\Delta Q-\Delta Q\delta^n Q-\Delta\delta^n QQ),\\
~~\mathrm{div}\delta^n u=0,\\
~~(\delta^n u,\delta^n Q)|_{t=0}=(u^n_0-u_0,Q^n_0-Q_0),~~(\delta^n u,\partial_\nu\delta^n Q)|_{\partial U}=(0,\mathcal{O}).
\end{cases}
\end{equation}
Now,
 in the same way as in the proof of  the Theorem 2.2 in \cite{MR3582195}, we multiply $-\Delta \delta^n Q+\delta^n Q$ with the first equation, then we take  the trace and integrate over the space $U$. Then
we multiply $\delta^n u$ to the second equation, integrate over the domain $U$ and we  add these two inequalities. Now by the same computations as in the proof  the Theorem 2.2 in \cite{MR3582195} (replacing $Q_1$ by $Q^n$, $Q_2$ by $Q$  $u_1$ by $u^n$ and $u_2$ by $u$), we obtain
\begin{equation*}
 \aligned
  \frac{1}{2}\frac{d}{dt}&\left(\|\delta^n Q\|^2_{L^2}+\|\nabla \delta^n Q\|^2_{L^2}+\|\delta^n u\|^2_{L^2}\right)+\|\Delta \delta^n Q\|^2_{L^2}+\|\nabla \delta^n u\|^2_{L^2}\\
  &+(a+1)\|\nabla \delta^n Q\|^2_{L^2}+a\|\delta^n Q\|^2_{L^2}\\
 \lesssim\,\,& \|\Delta\delta^n Q\|_{L^2}
\|\delta^n u\|_{L^2}\|\nabla Q\|_{L^\infty}+   \|\Delta\delta^n Q\|_{L^2}     \|u\|_{L^\infty}\|\nabla \delta^n Q\|_{L^2}\\
  &+  \|\delta^n u\|_{L^2}  \|\delta^n Q\|_{L^2}\|\nabla Q\|_{L^\infty}+
\|\Delta\delta^n Q\|_{L^2}\|\delta^n Q\|_{L^2}\|\nabla u\|_{L^\infty}+\|\nabla\delta^n u\|_{L^2}\|Q\|_{L^\infty}\|\delta^n Q\|_{L^2}\\
  &+|b|\|\Delta\delta^n Q\|_{L^2}\|\delta^nQ\|_{L^4} \|Q^n\|_{L^4}+|b|\|\Delta\delta^n Q\|_{L^2}\|\delta^n Q\|_{L^2}\|Q\|_{L^\infty}+|b| \|Q^n\|_{L^2}\|\delta^n Q\|^2_{L^4}\\
  &+|b| \|Q\|_{L^2}\|\delta^n Q\|^2_{L^4}
+c\|\Delta\delta^n Q\|_{L^2}\|\delta^n Q\|_{L^4}\|Q^n\|^2_{L^8}  \\
&   +c\|\Delta\delta^n Q\|_{L^2}\|\delta^n Q\|_{L^4}\|Q\|_{L^\infty}(\| Q_n\|_{L^4}+\|Q\|_{L^4})\\
&+c\| Q\|_{L^\infty}\|\delta^n Q\|^2_{L^4}(\| Q_n\|_{L^2}+\|Q\|_{L^2})\\
&+\|\nabla\delta^n Q\|_{L^2}\|\nabla Q\|_{L^\infty}\|\nabla \delta^n u\|_{L^2}+\|\Delta Q\|_{L^\infty}\|\delta^n Q\|_{L^2}\|\nabla \delta^n u\|_{L^2}\\ \leq&
\,\,\frac12\|\Delta \delta^n Q\|^2_{L^2}+\frac12\|\nabla \delta^n u\|^2_{L^2}+\frac a2\|\delta^n Q\|^2_{L^2}\\
&+C\| \delta^n u\|^2_{L^2}\underbrace{\left(\|\nabla Q\|^2_{L^\infty}+\|\nabla u\|^2_{L^\infty}\right)}_{J_1}+C\| \nabla\delta^n Q\|^2_{L^2}\underbrace{\left(\| u\|^2_{L^\infty}+\|\nabla Q\|^2_{L^\infty}\right)}_{J_2}\\
&+C\| \delta^n Q\|^2_{L^2}\underbrace{\left(1+\|\nabla u\|^2_{L^\infty}+\|Q\|^2_{L^\infty}+\|\Delta Q \|^2_{L^\infty}\right)}_{J_3}\\
&+C\| \delta^n Q\|^2_{L^4}\underbrace{\left(\|Q^n\|^2_{L^4}+\|Q^n\|_{L^2}+\|Q\|_{L^2} +\|Q^n\|^4_{L^8}\right)}_{J_4}\\
&+C\| \delta^n Q\|^2_{L^4}\underbrace{\|Q\|_{L^\infty}\left(\|Q^n\|_{L^2}+\|Q\|_{L^2}+\|Q\|_{L^\infty}(\|Q^n\|^2_{L^4} +
	\|Q\|^2_{L^4})\right)}_{J_5}.\\
\endaligned
\end{equation*}
Let
\begin{equation*}
  E_n(t):=\|\delta^n Q\|^2_{L^2}+\|\Delta \delta^n Q\|^2_{L^2}+\|\delta^n u\|^2_{L^2},
\end{equation*}
\begin{equation*}
  E_n(0):=\|Q^n_0-Q_0\|^2_{L^2}+\|\Delta Q^n_0-\Delta Q_0\|^2_{L^2}
+\|u^n_0-u_0\|^2_{L^2},
\end{equation*}
\begin{equation*}
\aligned
A_n(t):=\sum_{i=1}^{5}J_i
\endaligned
\end{equation*}
Therefore,
since
\begin{equation*}
\|\delta^n Q\|^2_{L^4}\leq C\|\delta^n Q\|^\frac{1}{2}_{L^2}\|\nabla\delta^n Q\|^\frac{3}{2}_{L^2}\leq C\left(\|\delta^n Q\|^{2}_{L^2}+\|\nabla\delta^n Q\|^{2}_{L^2}\right),
\end{equation*}
we have
\begin{equation*}
  \frac{d}{dt}E_n(t)\leq CA_n(t)E_n(t),
\end{equation*}
where $C$ is indipendent of time.
Since  ${H}_n(t),$ ${H}(t)$
are uniformly bounded on $[0,T]$ and using the embedding relation $W^{1,q}\hookrightarrow L^\infty,~~q>3$,  we have that $A_n(t)$
is integrable on $[0,T],$  and $\int^T_0A_n(s)ds$ is bounded by a constant that is indipendent of $n$. It follows from the Gronwall inequality that for any $t\in (0,T]$
\begin{equation*}
  E_n(t)\leq CE_n(0)\exp(C\int^t_0A_n(s)ds)\leq CE_n(0).
\end{equation*}
Since $E_n(0)\rightarrow0,$ we get $E_n(t)\rightarrow0$, as $ n\to \infty$ uniformely on $[0,T]$.
This ends the proof of {\sc Theorem b}.
\end{proof}

\par\noindent
{\bf Proof of {\sc Theorem c}.}
\par
\begin{proof}
Firstly, we introduce some dimensionless quantities:
\begin{equation*}
\aligned
& E_p(u,r):=\frac{1}{r^{5-p}}\int_{Q_r(z_0)}|u|^pdxds,~~~~E_{\ast}(u,r)=\frac{1}{r}\int_{Q_r(z_0)}|\nabla u|^2dxds,\\
& E(u,r):=\sup_{t_0-r^2\leq t<t_0}\frac{1}{r}\int_{B_r(x_0)}|u|^2dx,~~~~P_{\frac{3}{2}}(P,r):=\frac{1}{r^2}\int_{Q_r(z_0)}|P|^{\frac{3}{2}}dxds,\\
&E_*(u,\nabla Q;r)=E_*(u,r)+E_*(\nabla Q,r),~~~E(u,\nabla Q;r):=E(u;r)+E(\nabla Q;r),\\
&E_p(u,\nabla Q;r)=E_p(u,r)+E_p(\nabla Q,r)
\endaligned
\end{equation*}
where $z_0=(x_0,t_0)\in\R^3\times(0,\infty)$ and $Q_r(z_0)=Q_r(x_0,t_0):=B_r(x_0)\times[t_0-r^2,t_0].$
We will prove that the following dimensionless quantities are small:
\begin{equation}\label{CWW2}
  E(u,\nabla Q;r)+E_{\ast}(u,\nabla Q;r)+P_{\frac{3}{2}}(P,r)\leq \epsilon^2_1,
\end{equation}
when $r$ is enough small. Since $E_{\ast}(u,\nabla Q;r)=\frac{1}{r}\int_{Q_r(z_0)}|\nabla u|^2dxds+\frac{1}{r}\int_{Q_r(z_0)}|\nabla^2 Q|^2dxds$ the inequality \eqref{CWW2} would imply \eqref{du2} hand so thanks to Theorem \ref{Th1000} the solution $(u,Q)$ is smooth near $z_0$.
Therefore, it is enough  to prove \eqref{CWW2}.  We first remarke that  in \cite{Seregin,MR3904059}, the following inequalities are established
\begin{equation}\label{CWW5}
  E_3(u,r)\leq C\|u\|^{\frac{3}{2}}_{L^\infty(t_0-r^2,t_0;\dot{B}^{-1}_{\infty,\infty})}\bigg(E^{\frac{3}{4}}(u,2r)+E^{\frac{3}{4}}_{\ast}(u,2r)\bigg),
\end{equation}
and
\begin{equation}\label{CWW6}
  E_3(\nabla  Q,r)\leq C\|\nabla  Q\|^{\frac{3}{2}}_{L^\infty(t_0-r^2,t_0;\dot{B}^{-1}_{\infty,\infty})}\bigg(E^{\frac{3}{4}}(\nabla  Q,2r)+E^{\frac{3}{4}}_{\ast}(\nabla  Q,2r)\bigg).
\end{equation}
Let $\varphi\in C^\infty_0(Q_{2r}(x_0,t_0))$ be a function such that
\begin{equation*}
  \varphi=1 ~~\text{in}~Q_{r}(x_0,t_0),~~|\nabla \varphi|\leq \frac{1}{2r},~~|\nabla^2\varphi|+|\varphi_t|\leq \left(\frac{1}{2r}\right)^2.
\end{equation*}
We choose the $\phi=\varphi^2$ in the local energy inequality \eqref{SU} and after some calculations, we get
\begin{equation}\label{cz44}
\aligned
\sup_{t_0-4r^2\leq t\leq t_0}&\int_{\R^3}(|u|^2+|\nabla Q|^2)\varphi^2dx+\int_{\R^3\times[t_0-4r^2,t_0]}(|\nabla u|^2+|\nabla^2 Q|^2)\varphi^2dxdt\\
&\leq\int_{\R^3\times[t_0-4r^2,t_0]}(|u|^2+|\nabla Q|^2)(|\varphi_t||\varphi|+|\nabla \varphi|^2+|\nabla^2\varphi||\varphi|)dxdt\\
&+\int_{\R^3\times[t_0-4r^2,t_0]}[|u|^2+|\nabla Q|^2+|P|]|u||\nabla \varphi||\varphi|dxdt\\
&+({\|Q(t,\cdot)\|^2_{L^\infty(\R^3\times[t_0-4r^2,t_0])}+\|Q(t,\cdot)\|_{L^\infty(\R^3\times[t_0-4r^2,t_0])}})\int_{\R^3\times[t_0-4r^2,t_0]}|\nabla Q|^2\varphi^2dxdt\\
&+\|Q\|_{L^\infty(\R^3\times[t_0-4r^2,t_0])}\int_{\R^3\times[t_0-4r^2,t_0]}(|\nabla u||\nabla Q|+|u||\Delta Q|)|\varphi||\nabla \varphi|dxdt.
\endaligned
\end{equation}
In particolar, in order to establish \eqref{cz44}, we use the following estimates
\begin{equation*}
\aligned
  \int_{\R^3\times[t_0-4r^2,t_0]}[Q,\Delta Q]&\cdot u\otimes\nabla\varphi^2(x,s) dxds
  \leq\|Q\|_{L^\infty(\R^3\times[t_0-4r^2,t_0])}\int_{\R^3\times[t_0-4r^2,t_0]} |u||\Delta Q||\varphi||\nabla \varphi|dxdt
\endaligned
\end{equation*}
and
\begin{equation*}
\aligned
  \int_{\R^3\times[t_0-4r^2,t_0]}\nabla\mathcal{L}[\partial F(Q)]\cdot\nabla Q\varphi^2 dxds&\leq(\|Q(t,\cdot)\|^2_{L^\infty(\R^3\times[t_0-4r^2,t_0])}\\
  &~~~~+\|Q(t,\cdot)\|_{L^\infty(\R^3\times[t_0-4r^2,t_0])})\int_{\R^3\times[t_0-4r^2,t_0]}|\nabla Q|^2\varphi^2dxdt.
\endaligned
\end{equation*}
We now go on estimating  the second term in the inequality \eqref{cz44}.
In the last line of \eqref{cz44}, with the help of Young's inequality,  we obtain
\begin{equation}\label{cz55}
\aligned
  \int_{\R^3\times[t_0-4r^2,t_0]}&(|\nabla u||\nabla Q|+|u||\Delta Q|)|\varphi||\nabla \varphi|dxdt\\
  &\leq\frac{1}{2\|Q\|_{L^\infty(\R^3\times[t_0-4r^2,t_0])}}\int_{\R^3\times[t_0-4r^2,t_0]}(|\nabla u|^2+|\nabla^2 Q|^2)\varphi^2dxdt\\
  &~~+\|Q\|_{L^\infty(\R^3\times[t_0-4r^2,t_0])}\int_{\R^3\times[t_0-4r^2,t_0]}(|u|^2+|\nabla Q|^2)|\nabla \varphi|^2dxdt,
\endaligned
\end{equation}
and so the first right hand side term of the inequality \eqref{cz55} will be absorbed by the left hand side term of the inequality \eqref{cz44}.

Divide by $r$  both sides of inequality \eqref{cz44} and choose  $r$ such that  $r^2\leq\frac{1}{8}$ in order to obtain
\begin{equation}\label{cz66}
  \frac{1}{r}\int_{\R^3\times[t_0-4r^2,t_0]}|\nabla Q|^2\varphi^2dxdt\leq\frac{1}{(2r)^3}\int_{B_{2r}(x_0,t_0)}|\nabla Q|^2dxdt.
\end{equation}
Then from the definition of $E(u,\nabla Q;r), E_{\ast}(u,\nabla Q;r)$ and inequalities \eqref{cz44}, \eqref{cz55}, \eqref{cz66}, we get
\begin{equation*}
\aligned
&E(u,\nabla Q;r)+E_{\ast}(u,\nabla Q;r)\\
&\leq\frac{1}{(2r)^3}({\|Q(t,\cdot)\|^2_{L^\infty(\R^3\times[t_0-4r^2,t_0])}}+\|Q(t,\cdot)\|_{L^\infty(\R^3\times[t_0-4r^2,t_0])}+1)\int_{Q_{2r}(x_0,t_0)}(|u|^2+|\nabla Q|^2)dxdt\\
&~~~~~~~+\frac{1}{(2r)^2}\int_{Q_{2r}(x_0,t_0)}|P||u|dxdt+\frac{1}{(2r)^2}\int_{Q_{2r}(x_0,t_0)}(|u|^2+|\nabla Q|^2)|u|dxdt\\
&:=I_1+I_2+I_3.
\endaligned
\end{equation*}
It is easy to check that
\begin{equation*}
  |I_1|\leq ({\|Q(t,\cdot)\|^2_{L^\infty(\R^3\times[t_0-4r^2,t_0])}+\|Q(t,\cdot)\|_{L^\infty(\R^3\times[t_0-4r^2,t_0])}}+1)E^{\frac{2}{3}}_3(u,\nabla Q;2r),
\end{equation*}
and
\begin{equation*}
  |I_2|\leq E_3(u,\nabla Q;2r)+P_{\frac{3}{2}}(P,2r).
\end{equation*}
We estimate the term $I_3$ as follows:
\begin{equation*}
\aligned
  |I_3|\leq &E_3(u,\nabla Q;2r)+\frac{1}{(2r)^2}\bigg(\int_{Q_{2r}(x_0,t_0)}|\nabla Q|^3dxdt\bigg)^{\frac{2}{3}}\bigg(\int_{Q_{2r}(x_0,t_0)}|u|^3dxdt\bigg)^{\frac{1}{3}}\\
  &\leq CE_3(u,\nabla Q;2r).
\endaligned
\end{equation*}
We remark that $\|Q(t,\cdot)\|_{L^\infty(\R^3\times[t_0-4r^2,t_0])}\leq C\|Q_0\|_{L^\infty}$ see [21], so thanks to  \eqref{CWW5} e \eqref{CWW6}, we have
\begin{equation}\label{CWW3}
\aligned
  &E(u,\nabla Q;r)+E_{\ast}(u,\nabla Q;r)\leq CE^{\frac{2}{3}}_3(u,\nabla Q;2r)+CE_3(u,\nabla Q;2r)+CP_{\frac{3}{2}}(P,2r)\\
  &\leq \|(u,\nabla Q)\|_{L^\infty(t_0-4r^2,t_0;\dot{B}^{-1}_{\infty,\infty})}\bigg(E^{\frac{1}{2}}(u,4r)+E^{\frac{1}{2}}_{\ast}(u,4r)\bigg)\\
  &~~~+\|(u,\nabla Q)\|^{\frac{3}{2}}_{L^\infty(t_0-4r^2,t_0;\dot{B}^{-1}_{\infty,\infty})}\bigg(E^{\frac{3}{4}}(u,4r)+E^{\frac{3}{4}}_{\ast}(u,4r)\bigg)\\
  &~~~+CP_{\frac{3}{2}}(p,2r).
\endaligned
\end{equation}
Then, we prove the following estimate of $P_{\frac{3}{2}}(P,r)$ similar to the Lemma 6.3 in paper \cite{MR4134151}.
\begin{equation}\label{CWW4}
  P_{\frac{3}{2}}(P,r)\leq C\left(\Big(\frac{\rho}{r}\Big)^2E_3(u,\nabla Q;\rho)+\Big(\frac{r}{\rho}\Big)P_{\frac{3}{2}}(P,\rho)\right)~~~\text{when}~0<r\leq\frac{\rho}{2}.
\end{equation}
From the scaling invariance of all quantities, we only need to consider the case $\rho=1,$ $0<r\leq\frac{1}{2}$, $z_0=(0,0).$
Taking the divergence of the second equation in \eqref{cz22} and using the fact [Lemma 2.3 in \cite{MR4134151}],
\begin{equation*}
  \mathrm{div}^2[Q,H]=0, ~~\text{in}~~ \R^3,
\end{equation*}
we obtain
\begin{equation*}
\aligned
-\Delta P=\mathrm{div}^2[u\otimes u+\nabla Q\otimes\nabla Q]
\endaligned
\end{equation*}
Let $\eta\in C^\infty_0(\R^3)$ be a cut-off function of $B_{\frac{2}{3}}(0)$ such that
\begin{equation*}
  \eta=1,~~\text{in}~B_{\frac{2}{3}}(0),~~\eta=0,~~\text{in}~~\R^3\setminus B_1(0),~~0\leq\eta\leq1,~~|\nabla\eta|\leq8.
\end{equation*}
Define the following auxiliary function:
\begin{equation*}
\aligned
P_1(x,t)=-\int_{\R^3}\nabla^2_yG(x-y):\eta^2(y)[u\otimes u+\nabla Q\otimes\nabla Q](y,t)dy.
\endaligned
\end{equation*}
where $G$ is the fundamental solution of harmonic equation.
It is easy to see that the $P_1$ is solution of the following equation:
\begin{equation*}
  -\Delta P_1=\mathrm{div}^2[u\otimes u+\nabla Q\otimes\nabla Q]~~\text{in}~B_{\frac{2}{3}}(0),
\end{equation*}
and
\begin{equation*}
  -\Delta(P-P_1)=0~~\text{in}~B_{\frac{2}{3}}(0).
\end{equation*}
We estimate the $P_1$ by the Calderon-Zygmund theory:
\begin{equation}\label{FF1}
\aligned
\|P_1\|^{\frac{3}{2}}_{L^{\frac{3}{2}}(\R^3)}&\leq \|\eta^2|u|^2\|^{\frac{3}{2}}_{L^{\frac{3}{2}}(\R^3)}+\|\eta^2|\nabla Q|^2\|^{\frac{3}{2}}_{L^{\frac{3}{2}}(\R^3)}\\
&\leq\|u\|^3_{L^3(B_1(0))}+\|\nabla Q\|^3_{L^3(B_1(0))}.
\endaligned
\end{equation}
Since $P-P_1$ is harmonic in $B_{\frac{2}{3}}(0),$ using the mean-value property, one has
\begin{equation}\label{FF2}
\aligned
\frac{1}{r^2}&\|P-P_1\|^{\frac{3}{2}}_{L^{\frac{3}{2}}(B_r(0))}=\frac{1}{r^2}\int_{B_r(0)}|P-P_1|^\frac{3}{2}dx\\
&\leq \frac{1}{r^2}r^3\|P-P_1\|^\frac{3}{2}_{L^{\infty}(B_{r}(0))}\\
&\leq \frac{1}{r^2}r^3\sup_{x\in B_{r}(0)}\bigg(\frac{1}{B_{\frac{1}{8}}(x)}\int_{B_{\frac{1}{8}}(x)}|P-P_1|dx\bigg)^\frac{3}{2}\\
&\leq r\int_{B_1(0)}|P-P_1|^\frac{3}{2}dx\\
&\leq r(\|P\|^{\frac{3}{2}}_{L^{\frac{3}{2}}(B_1(0))}+\|P_1\|^{\frac{3}{2}}_{L^{\frac{3}{2}}(B_1(0))}).
\endaligned
\end{equation}
where we have used the fact $B_{\frac{1}{8}}(x)\subset B_\frac{2}{3}(0),$ for $x\in B_{r}(0)$ since $r\leq\frac{1}{2}.$
Integrating $\|P\|^{\frac{3}{2}}_{L^{\frac{3}{2}}(B_1(0))}$ over $[-r^2,0]$ and using the previous estimates \eqref{FF1}, \eqref{FF2}, we can show
\begin{equation*}
\aligned
\frac{1}{r^2}\int_{Q_r(0,0)}|P|^{\frac{3}{2}}dxds\leq CrP_{\frac{3}{2}}(P,1)+C\frac{1}{r^2}E_3(u,\nabla Q;1),~~\text{with}~~\rho=1,
\endaligned
\end{equation*}
so we complete the proof of the inequality \eqref{CWW4}.
It follows from \eqref{CWW5}, \eqref{CWW6} and \eqref{CWW4} that
\begin{equation}\label{FF4}
\aligned
P_{\frac{3}{2}}(P,r)&\leq C\left(\Big(\frac{\rho}{r}\Big)^2E_3(u,\nabla Q;\rho)+\frac{r}{\rho}P_{\frac{3}{2}}(P,\rho)\right)\\
&\leq C\left(\Big(\frac{\rho}{r}\Big)^2\|(u,\nabla Q)\|^{\frac{3}{2}}_{L^\infty(t_0-\rho^2,t_0;\dot{B}^{-1}_{\infty,\infty})}\Big(E^{\frac{3}{4}}(u,2\rho)+E^{\frac{3}{4}}_{\ast}(u,2\rho)\Big)
+\frac{r}{\rho}P_{\frac{3}{2}}(P;2\rho)\right)\\
&\leq C\left(\Big(\frac{\varrho}{r}\Big)^2\|(u,\nabla Q)\|^{\frac{3}{2}}_{L^\infty(t_0-\rho^2,t_0;\dot{B}^{-1}_{\infty,\infty})}\Big(E^{\frac{3}{4}}(u,\varrho)+E^{\frac{3}{4}}_{\ast}(u,\varrho)\Big)
+\frac{r}{\varrho}P_{\frac{3}{2}}(P;\varrho)\right),
\endaligned
\end{equation}
where $\varrho:=2\rho.$
\par\noindent
Set $F_1(r):=E(u,\nabla Q,r)+E_{\ast}(u,\nabla Q;r)+P_{\frac{3}{2}}(P,r)$ and  $M(r):=\|(u,\nabla Q)\|_{L^\infty(t_0-r^2,t_0;\dot{B}^{-1}_{\infty,\infty})}$, $M:=M(\tilde{\varrho})$.
Let $\tilde{\varrho}$ be a fixed positive constant, then for $r$ small enough, we have
\begin{equation*}
  E^{\frac{1}{2}}(u,4r)\leq(\frac{\tilde{\varrho}}{r})^{\frac{1}{2}}E^{\frac{1}{2}}(u,\tilde{\varrho}),
\end{equation*}
and so \eqref{CWW3} and \eqref{FF4} with $\varrho=\tilde{\varrho}$ implies that
for fixed $\tilde{\varrho}$ and enough small $r$
 \begin{equation}\label{cw100}
\aligned
F_1(r)&\leq C\Bigg( M(\frac{\tilde{\varrho}}{r})^{\frac{1}{2}}(E^{\frac{1}{2}}(u,\nabla Q;\tilde{\varrho})+E^{\frac{1}{2}}_*(u,\nabla Q;\tilde{\varrho}))\\
  &+M^\frac{3}{2}(\frac{\tilde{{\varrho}}}{r})^{\frac{3}{4}}(E^{\frac{3}{4}}(u,\nabla Q;\tilde{\varrho})+E^{\frac{3}{4}}_*(u,\nabla Q;\tilde{\varrho}))\\
  &+M^\frac{3}{2}(\frac{\tilde{\varrho}}{r})^2\bigg(E^{\frac{3}{4}}(u,\tilde{\varrho})+E^{\frac{3}{4}}_{\ast}(u,\tilde{\varrho})\bigg)
+\frac{r}{\tilde{\varrho}}P_{\frac{3}{2}}(P;\tilde{\varrho})\Bigg)\\
&\leq C\Bigg( (M+M^\frac{3}{2})\Big(\frac{\tilde{\varrho}}{r}\Big)^2(F_1^{\frac{3}{4}}(\tilde{\varrho})+
F_1^{\frac{1}{2}}(\tilde{\varrho}))+\frac{r}{\tilde{\varrho}}
(F_1(\tilde{\varrho})
)\Bigg)\\
&\leq C\Bigg(\frac{r}{\tilde{\varrho}} F_1(\tilde{\varrho})+(M+M^\frac{3}{2})^4\Big(\frac{\tilde{\varrho}}{r}\Big)^{11}\Bigg).
\endaligned
\end{equation}
where in the last inequality we control the term $(M+M^\frac{3}{2})(\frac{\tilde{\varrho}}{r})^2(F_1^{\frac{3}{4}}(\tilde{\varrho})+F_1^{\frac{1}{2}}(\tilde{\varrho}))$  by
$ C(\frac{r}{\tilde{\varrho}})F_1(\tilde{\varrho})+C(M+M^\frac{3}{2})^4(\frac{\tilde{\varrho}}{r})^{11}$ by Young inequality.
Note the constant $C$ in equality \eqref{CWW4} is independent on the ratio $r$ with $\rho$. Set $\theta:=\frac{r}{\tilde{\varrho}}$ with $r$ small enough such that
$C\theta<1.$
Then from the inequality \eqref{cw100} it follows that
\begin{equation*}
  F_1(\theta\tilde{\varrho})\leq C\theta F_1(\tilde{\varrho})+C\frac{(M+M^\frac{3}{2})^4}{\theta^{11}},
\end{equation*}
Iterating the above inequality, we get
\begin{equation}\label{EC12}
  F_1(\theta^k\tilde{\varrho})\leq (C\theta)^k
  F_1(\tilde{\varrho})+C(M+M^\frac{3}{2})^4(\frac{1}{\theta})^{11}\frac{1-(C\theta)^k}{1-(C\theta)}
\end{equation}
Now we prove that there exists $k_0\in N$ such that $\forall k\geq k_0,$
\begin{equation}\label{YY}
 F_1(\theta^k\tilde{\varrho})\leq \epsilon^2_1\theta^2,
\end{equation}
where $\epsilon_1$ is the constant appearing in the Theorem \ref{Th1000}.
We remark that the following embedding holds
\begin{equation*}
  \|\nabla Q\|_{L^\infty(t_0-r^2,t_0;\dot{B}^{-1}_{\infty,\infty})}\leq\|Q\|_{L^\infty(t_0-r^2,t_0;\dot{B}^{0}_{\infty,\infty})}\leq
  \|Q(t)\|_{L^\infty(t_0-r^2,t_0;L^\infty)}\leq e^{ct}\|Q_0\|_{L^\infty},
\end{equation*}
thanks to the property
\begin{equation*}
  \|Q(t,\cdot)\|_{L^\infty}\leq e^{ct}\|Q_0\|_{L^\infty},~~~\forall~t\geq 0.
\end{equation*}
 proved in \cite{MR2864407}. So using the hypothesis \eqref{cww1} and \eqref{YYY2},
we get
\begin{equation*}
  M\leq \varepsilon_0+e^{ct_0}\|Q_0\|_{L^\infty}\leq\varepsilon_0+e^{ct_0}\varepsilon_1.
\end{equation*}
Therefore, choosing $k_0$ large enough, and $\varepsilon_0, \varepsilon_1$ small enough, we get \eqref{YY}.
 For any $0<r<\theta^{k_0}\tilde{\rho},$
we can find $k_1$ such that $\theta^{k_1+1}\tilde{\varrho}\leq r\leq\theta^{k_1}\tilde{\varrho},$ From the definition of $F_1,$ one gets
\begin{equation*}
  F_1(r)\leq \frac{1}{\theta^2}F_1(\theta^{k_1}\tilde{\varrho})\leq \epsilon^2_1.
\end{equation*}
Moreover, we can get the solution $(u,Q)$ is smooth near the point $(x_0,t_0)\in\R^3\times(0,+\infty)$ by Theorem \ref{Th1000}.
So we complete the proof of {\sc Theorem c}.
\end{proof}
\vskip4mm

\section*{Acknowledgments}
This work is partially supported by the National Natural Science Foundation of China (Nos. 11801574, 11971485), Natural Science Foundation of Hunan Province (No. 2019JJ50788), Central South University Innovation-Driven Project for Young Scholars (No. 2019CX022) and Fundamental Research Funds for the Central Universities of Central South University, China (Nos. 2020zzts038, 2021zzts0041).

\par\noindent
\bibliographystyle{siam}
\bibliography{Reference}
\end{document}